\documentclass[12pt,fleqn]{article}

\usepackage{mathrsfs}
\usepackage{graphicx}
\usepackage{amsmath}
\usepackage{amsfonts}
\usepackage{amssymb}

\newtheorem{theorem}{Theorem}[section]
\newtheorem{lemma}[theorem]{Lemma}

\numberwithin{equation}{section}

\begin{document}

\title{ Tricyclic graphs with exactly two main eigenvalues \thanks{This research
was partially supported by the NSF of China(No.10971086).}}

\author{Xiaoxia Fan, Yanfeng Luo\\
{\small Department of Mathematics, Lanzhou University,}\\
{\small Lanzhou, Gansu 730000, PR China}\\
{\small fanxx06@lzu.cn}}
\date{}

\maketitle

\begin{abstract}
An eigenvalue of a graph $G$ is called a main eigenvalue if it has
an eigenvector the sum of whose entries is not equal to zero. In
this paper, all connected tricyclic graphs with exactly two main
eigenvalues are determined.

\vskip 0.05in


\end{abstract}

\section{ Introduction }
All graphs considered in this paper are finite, undirected and
simple. Let $G = (V,E)$ be a graph with vertex set $V(G)$ and edge
set $E(G)$. Denote by $A(G)$  the adjacency matrix of $G$. The
eigenvalues of $G$ are those of $A(G)$. An eigenvalue of a graph $G$
is called a main eigenvalue if it has an eigenvector the sum of
whose entries is not equal to zero. It is well known that a graph is
regular if and only if it has exactly one main eigenvalue.

A long-standing problem posed by Cvetkovic (\cite{Cvetkovic})
is that of how to characterize graphs with exactly $k (k\geq 2)$
main eigenvalues. Hagos \cite{Hagos} gave an alternative
characterization of graphs with exactly two main eigenvalues.
Recently, Hou and Zhou  \cite{Hou} characterized the trees with
exactly two main eigenvalues.

A vertex of a graph $G$ is said to be pendant if it has degree one. Denote by
$C_n$ and $P_n$ the cycle
and path of order $n$, respectively. A connected graph is said to be tricyclic
(resp., unicyclic and bicyclic), if  $|E(G)|=|V(G)|+2$ (resp., $|E(G)|=|V(G)|$
and $|E(G)|=|V(G)|+1$). Hou and Tian
\cite{YHou} showed that the graphs $C_{r}^{k}$ for some positive integers $k, r$
with $r\geq 3$, where $C_{r}^{k}$ is the graph obtained from $C_r$ by attaching
$k > 0$ pendant vertices to every
vertex of $C_r$, are the only connected unicyclic graphs with exactly
two main eigenvalues. Hu et al. \cite{Hu} and Shi \cite{Shi} characterized all
connected bicyclic graphs with exactly two main eigenvalues
independently. This paper will continue the line of this research and determine
all connected tricyclic graphs with exactly two main eigenvalues.

For any tricyclic graph $G$, the base of $G$,
denoted by $G_B$ is the minimal tricyclic subgraph of $G$. Clearly,
$G_B$ is the unique tricyclic subgraph of $G$ containing no pendent
vertex, and $G$ can be obtained  from  $G_B$
 by attaching trees to some vertices of $G_B$. It follows from \cite{Geng}
 that there are 8 types of bases for tricyclic graphs, say, $
{\mathcal T}_i,i=1,\dots,8$, which are depicted in Fig.~1.
\begin{figure}[!htbp]
\centering
\includegraphics[totalheight=4cm]{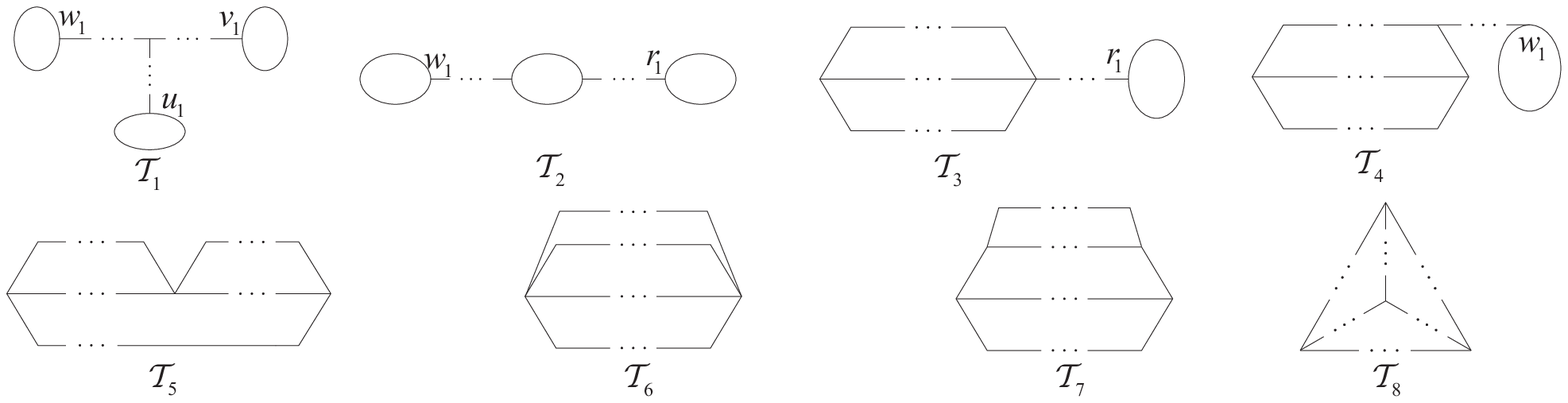}\\
 \caption{\label{lem} \small{The 8 types of bases for tricyclic graphs. }}
\end{figure}

\section{Preliminaries }
In this section, we will present some notations and known results
which will be used in the next section. The reader is referred to
\cite{Bondy} for any undefined notation and terminology on graphs in
this paper.

Let $G$ be a graph. As usual, we denote by $d(v)=d_G(v)$ and
$N(v)=N_G(v)$ the degree of vertex $v$ and the set of all neighbors
of $v$ in $G$. Let
\begin{equation}\label{s21}
S(v) = \sum\limits_{u\in N(v)}d(u).
\end{equation}
A graph $G$ is called 2-walk $(a,b)$-linear if there exist unique
rational numbers $a,b$ such that
\begin{equation}\label{sab}
S (v) = ad(v)+b
\end{equation}
holds for every vertex $v\in V(G)$.

An internal path of $G$ is a walk $v_0 v_1 \dots v_s$ such that the
vertices $v_0, v_1 ,\dots,v_s$ are distinct, $d(v_0 )>2,d(v_s )>2$,
and $d (v _i)= 2$ for $0<i<s$. An internal path is called an
internal cycle if $v_0 = v_s$. If $R$ is a path or a cycle of $G$,
the length of $R$, denoted by $l(R) $, is defined as the number of
edges of $R$.

\begin{lemma}(\cite{Hagos}). \label{222ab}
A graph $G$ has exactly two main eigenvalues if and only if $G$ is
2-walk $(a,b)$-linear.
\end{lemma}

\begin{lemma}(\cite{YHou}). \label{21} Let $G$ be a 2-walk
(a,b)-linear graph. Then both $a$ and $b$ are integers.
\end{lemma}

\begin{lemma}(\cite{Hu})\label{22}. Let $G$ be a 2-walk
$(a,b)$-linear graph and $v,u$ be two vertices of $G$ with unequal
degree $d(v),d(u)$, respectively. Then
\begin{equation} \label{suv}
a=\frac{S(v)-S(u)}{d(v)-d(u)}, \
b=\frac{d(u)S(v)-d(v)S(u)}{d(v)-d(u)}.
\end{equation}
\end{lemma}

In the following, for convenience, we always assume that ${\mathcal
G}$ is the set of tricyclic graphs with exactly two main
eigenvalues, $x$ is a pendent vertex of $G$ (if exist). For each
$G\in{\mathcal G}$, let $G_0$ be the graph obtained from $G$ by
deleting all pendent vertices. From the proof Lemmas 3.1-3.7 in
\cite{Hu}, we know that those Lemmas are also hold for  tricyclic
graphs. Hence we have the following two Lemmas.

\begin{lemma}\label{31}
Let $G\in {\mathcal G}$ and  $R=x_1x_2\dots x_t$ be an internal path
of length at least $2$ in $G$. Then $l(R)\leq3$. In particular, if $
l(R)=3$, then there exists no  path $Q=y_1y_2y_3$ in $G$ such that
$d(y_1)=d(y_3)=d(x_1)$ and $d(y_2)=2$.
\end{lemma}

\begin{lemma}\label{3no}
Let $G\in{\mathcal G}$ and $v\in V(G_{0})$.  Then

(i) $G_0\in {\mathcal T}_i,i=1,\dots,8$ (see Fig.~1);

(ii)  $d(v)=d_{G_0}(v)$ or $a+b$;

(iii) if $G$ has at least one pendent vertex, then $S(x)=a+b\geq3$
and
 $a\geq2$;

(iv) for a cycle $C=x_1x_2\dots x_tx_1$ of $G$ with
$d_{G_0}(x_1)\geq 3$, $d_{G_0}(x_2)=2$, if $G$ has at least one
pendent vertex, then there is an integer $i\in\{1,2,\dots,t\}$ such
that $d(x_i)\neq a+b$.
\end{lemma}

\section{Tricyclic graphs with exactly two main eigenvalues}

In this section, we will  determine all tricyclic graphs with
exactly two main eigenvalues. By Lemma \ref{222ab}, it is sufficient
to determine all 2-walk $(a,b)$-linear tricyclic graphs.

\begin{lemma}\label{p}
Let $G\in{\mathcal G}$ has at least one pendent vertex and let
$R=x_1x_2\dots x_t$ be an internal path or an internal cycle of
length at least $3$ in $G_0$ with $d_{G_0}(x_1)=d_{G_0}(x_t)\in
\{3,4,6\}$ or $d_{G_0}(x_1)=3,d_{G_0}(x_t)=5$. Then

(i) $ d(x_2)= d(x_3)=\dots=d(x_{t-1})\in\{2,a+b\}$ and
$d(x_1)=d(x_t)$;

(ii) if $d(x_2)=2$, then $l(R)=3$;

(iii) if $R$ is a cycle with $ d_{G_0}(x_1)\in\{3,4,6\}$, then
$l(R)=3$ and $d(x_2)=d(x_3)=2$. In particular, if $d_{G_0}(x_1)=3$,
then $a=2$;

(iv) if $R$ is a cycle with $ d_{G_0}(x_1)=5$, then $l(R)=3$,
$d(x_2)=d(x_3)\in\{2,3\}$. In particular, if $d(x_2)=d(x_3)=3$, then
$ d(x_1)=5 $ and $a=3,b=0$.
\end{lemma}
\noindent{\bf Proof.} (i) By way of contradiction, assume that there
is an integer $i\in\{2,3,\dots,t-2\}$ such that $d(x_i)\neq
d(x_{i+1})$. Without loss of generality, suppose that $i$ is the
smallest integer such that $d(x_i)\neq d(x_{i+1})$. By Lemma
\ref{3no}~(ii), we may assume that $d(x_i)=2$ and $d(x_{i+1})=a+b$.
Hence $d(x_2)=d(x_3)=\dots=d(x_i)=2$. Applying \eqref{suv} with
$(v,u)=(x_{i+1},x)$, we have
\begin{equation*}
a=\frac{S(x_{i+1})-S(x)}{d(x_{i+1})-d(x)}=\frac
{a+b-2+2+d(x_{i+2})-(a+b))}{a+b-1}=\frac{d(x_{i+2})}{a+b-1}.
\end{equation*}
This together with Lemma \ref{3no}~(iii) implies that
\begin{equation}\label{dxi2}
d(x_{i+2})=a(a+b-1)\geq2(a+b-1)\geq a+b+1>max\{a+b,3\}.
\end{equation}
By Lemma \ref{3no}~(ii), we have  $d(x_j)\in\{a+b,2\}$ for $2\leq
j\leq t-1$. Thus $x_{i+2}\in\{x_1,x_t\}$.

If $d_{G_0}(x_1)=d_{G_0}(x_t)=3$, then $d(x_{i+2})\in
\{d(x_1),d(x_t)\}\subseteq\{3,a+b\}$, contrary to \eqref{dxi2}.

If $d_{G_0}(x_1)=d_{G_0}(x_t)=4$, then $d(x_1),d(x_t)\in\{4,a+b\}$.
It follows from \eqref{dxi2} that $d(x_{i+2})=a(a+b-1)=4$. This
together with Lemma \ref{3no}~(iii) implies that $a=2, b=1$. Note
that $d(x_2)=2$. Then $S(x_2)=5$ by \eqref{sab}. On the other hand,
$d_{G_0}(x_1)=4>a+b$, so $d(x_1)=4$ by Lemma \ref{3no}~(ii). Thus by
\eqref{s21}, $S(x_2)=d(x_1)+d(x_3)=4+d(x_3)>5$, a contradiction.

If $d_{G_0}(x_1)=d_{G_0}(x_t)=6$, with a similar argument of the
case $d_{G_0}(x_1)=d_{G_0}(x_t)=4$, we will get a contradiction
again.

If $d_{G_0}(x_1)=3,d_{G_0}(x_t)=5$, then by \eqref{dxi2}
$d(x_{i+2})=d(x_t)=d_{G_0}(x_t)=5$ and $a(a+b-1)=5$. It contradicts
the fact that $a\geq 2, a+b\geq 3$.

Hence $d(x_2)= d(x_3)=\dots=d(x_{t-1})\in\{2,a+b\}$. Therefore
$S(x_2)=S(x_{t-1})$ by \eqref{sab}. On the other hand, by
\eqref{s21}, $S(x_2)=d(x_2)-2+d(x_1)+d(x_3),
S(x_{t-1})=d(x_{t-1})-2+d(x_t)+d(x_{t-3})$. It follows that
$d(x_1)=d(x_t)$.

(ii) Suppose contrary that $l(R)\geq4$. Then $d(x_3)=d(x_4)=2$ by
(i). Applying \eqref{suv} with $(v,u)=(x_{3},x)$, we have
$a=4-(a+b)$, which is impossible by Lemma \ref{3no}~(iii).

(iii) By Lemma \ref{3no}~(ii), we have $d(x_2)\in\{2,a+b\}$. If
$d(x_2)=a+b$, then $d(x_i)=a+b$ for $2\leq i\leq t-1$ by (i). Thus
$d(x_1)=d_{G_0}(x_1)\neq a+b$ by Lemma \ref{3no}~(iv). Applying
\eqref{suv} with $(v,u)=(x_2,x)$, we have
\begin{equation}\label{lcc}
a=\frac{a+b-2+d(x_1)+a+b-(a+b)}{a+b-1}=1+\frac{d(x_1)-1}{a+b-1}=1+\frac{d_{G_0}(x_1)-1}{a+b-1}.
\end{equation}
Note that $d_{G_0}(x_1)=3,4$ or $6$ and $d_{G_0}(x_1)\neq a+b\geq3$.
It follows from \eqref{lcc} that $a$ can not be an integer. It
contradicts Lemma \ref{21}. Hence $d(x_2)=2$. Therefore $l(R)=3$ and
$d(x_3)=d(x_2)=2$.

In particular, if $d_{G_0}(x_1)=3$, then $a=2+d(x_1)-(a+b)$ by
applying \eqref{suv} with $(v,u)=(x_2,x)$. If $d(x_1)=a+b$, then
$a=2$. If $d(x_1)=d_{G_0}(x_1)=3$, then $a=5-(a+b)$. It follows from
Lemmas \ref{21} and \ref{3no}~(iii) that $a=2$.

(iv) By Lemma \ref{3no}~(ii), $d(x_2)\in\{2,a+b\}$. If $d(x_2)=2$,
then $l(R)=3$ and $d(x_3)=2$ by~(i) and~(ii). If $d(x_2)=a+b$, then
$d(x_i)=a+b$ for $2\leq i\leq t-1$  by~(i). So
$d(x_1)=d_{G_0}(x_1)=5\neq a+b$ by Lemma \ref{3no}~(iv). Applying
\eqref{suv} with $(v,u)=(x_2,x)$, we have $a=1+\frac{4}{a+b-1}$.
Thus $a+b=3$ and $a=3$ by Lemmas \ref{21} and \ref{3no}~(iii).
Suppose that $l(R)\geq4$. Then $d(x_2)=d(x_3)=d(x_4)=a+b=3$ by (ii).
Applying \eqref{suv} with $(v,u)=(x_3,x)$, we have $a=2$. It is a
contradiction. Therefore $l(R)=3$. $\square$
\begin{figure}[!htbp]
\begin{center}
\includegraphics[totalheight=10.7cm]{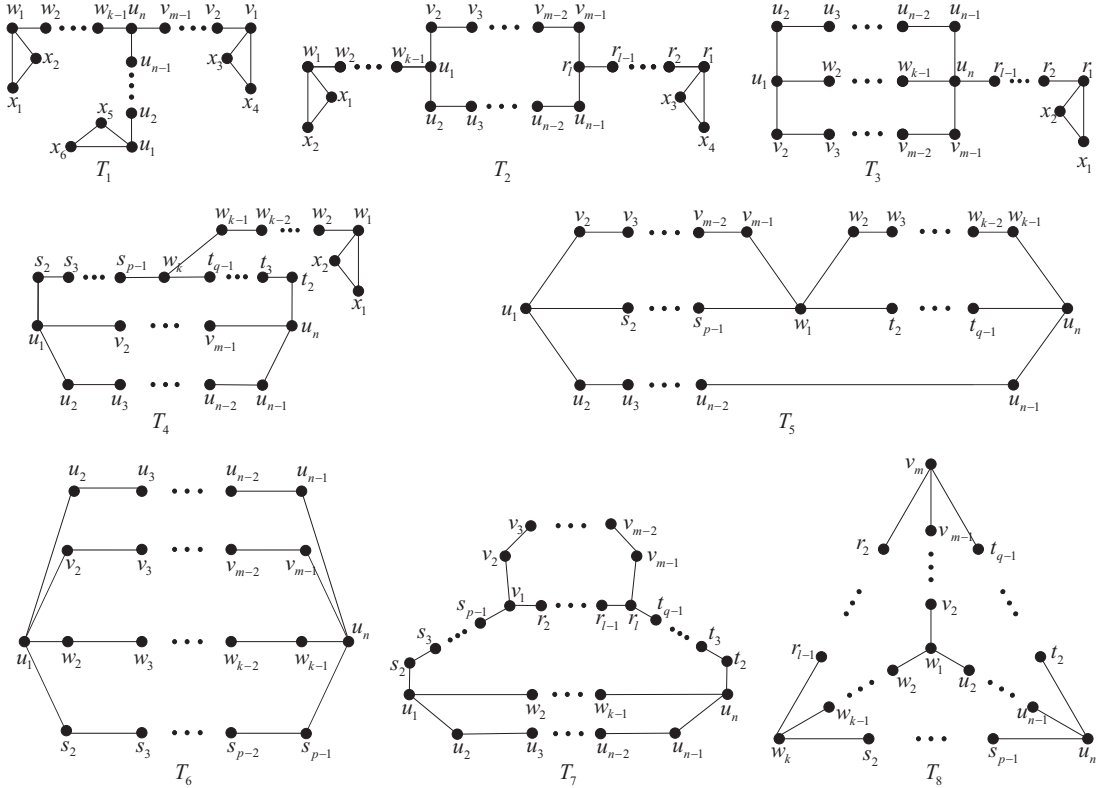}
 \caption{\label{lem} \small{The graphs $T_i$ for $i=1,\dots,8$.  }}
\end{center}
\end{figure}
\begin{figure}[!htbp]
\begin{center}
\includegraphics[totalheight=8cm]{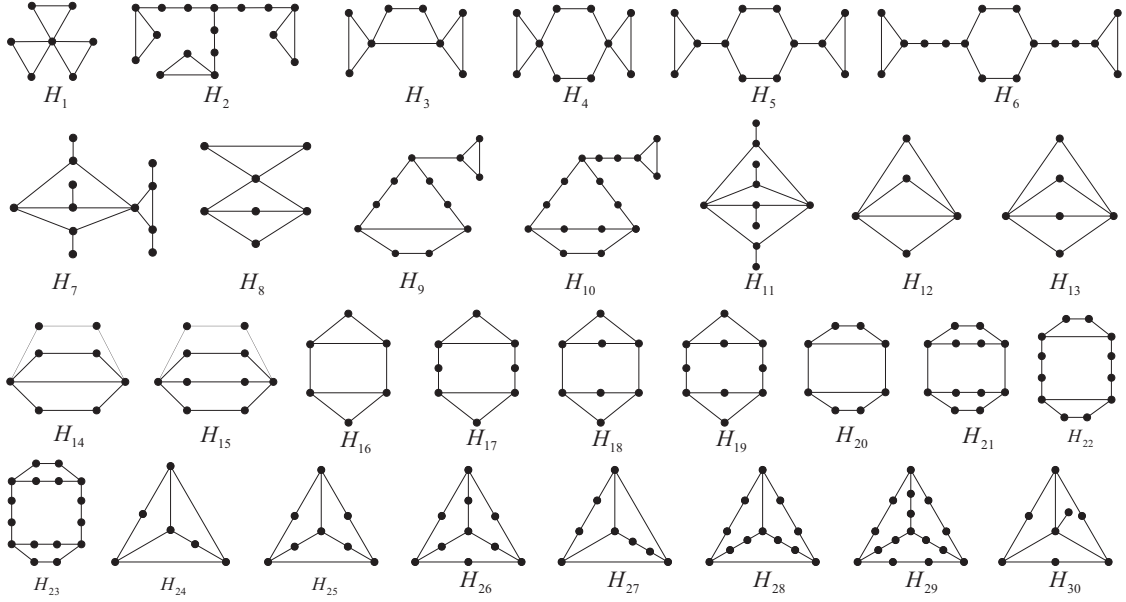}\\
 \caption{\label{main1} \small{The graphs $H_i$ for i=1,\dots,30. }}
\end{center}
\end{figure}
\begin{lemma}\label{t1} Let $G\in{\mathcal
G}$ with $G_0\in {\mathcal T}_1$ (see Fig.~1). Then $G=H_i$ for
$i=1,2$ (see Fig.~3).
\end{lemma}
\noindent {\bf Proof.} If $G_0\in {\mathcal T}_1$, then
$d_{G_0}(u_1), d_{G_0}(v_1), d_{G_0}(w_1)\in \{3, 4, 5, 6\}$ and so
each cycle of $G_0$ has the length of $3$ by Lemmas \ref{31} and
\ref{p}. Hence $G_0=T_1$ (see Fig.~2), where $n, m,k\geq 1$. For
convenience, we set $w_k=v_m=u_n$.

{\bf Case 1}. $n=m=k=1$.  If $G$ has no pendent vertex, then $G=
H_1$ (see Fig.~3). By \eqref{suv}, $H_1$ is 2-walk $(1,6)$-linear.
If $G$ has at least one pendent vertex, say $x$, then $x\in N(u_1)$
since $d(x_i)=2$ for $i=1,\dots,6$ by Lemma \ref{p}~(iii). It
follows from Lemma \ref{3no}~(ii) that $d(u_1)=a+b>6$. Applying
\eqref{suv} with $(v,u)=(u_1,x)$, we have
$a=\frac{a+b-6+12-(a+b)}{a+b-1}<2$. This contradicts Lemma
\ref{3no}~(iii).

{\bf Case 2}. $n\geq2$. We consider the following two cases:

{\bf Subcase 1}. $G$ has no pendent vertex. Then
$S(x_1)=2+d(w_1)=S(x_5)=5$ by \eqref{s21} and \eqref{sab}. Hence
$d(w_1)=3$. It implies that $k\geq2$. Similarly, we have $d(v_1)=3$
and $m\geq2$. By Lemma \ref{31}, we have $n,m,k\in\{2,4\}$. Without
loss of generality, suppose that $ n\geq m \geq k$. If $n=2$, then
$m=k=2$. Hence $S(u_1)=7,S(u_2)=9$  by \eqref{s21}. On the other
hand, $d(u_1)=d(u_2)=3$, so $S(u_1)=S(u_2)$ by \eqref{sab}, a
contradiction. Hence $n=4$. Similarly, we have $m=k=4$. Therefore
$G=H_2$ (see Fig.~3). By \eqref{suv}, $H_2$ is 2-walk
$(1,3)$-linear.

{\bf Subcase 2}. $G$ has at least one pendent vertex. In this case,
we show that there is no such graph with exactly two main
eigenvalues. Since $d_{G_0}(u_1)=3$, we have $a=2$ by Lemma
\ref{p}~(iii). So $d(x_i)=2$ for $i=1,\dots,6$ by Lemma
\ref{p}~(iii) and (iv).

We claim that $m, k\geq 2$. Otherwise, let $k=1$. Then
$S(x_1)=2+d(w_1)=S(x_5)=2+d(u_1)$ by \eqref{s21} and \eqref{sab}. It
follows from Lemma \ref{3no}~(ii) and the fact that
$d_{G_0}(u_1)\neq d_{G_0}(w_1)$ that $d(w_1)=d(u_1)=a+b$. Hence
$S(u_1)=S(w_1)$ by \eqref{sab}. By \eqref{s21},
\[
S(u_1)=a+b-3+4+d(u_2), S(w_1)=\Big\{\begin{array}{ll}
a+b-5+8+d(u_{n-1}), & \mbox{if}~ m=1,\\
a+b-4+4+d(v_{m-1})+d(u_{n-1}), &  \mbox{if}~ m\geq2.\\
\end{array}
\]
If $m=1$, then $d(u_2)=d(u_{n-1})+2$. This is impossible by Lemma
\ref{p}~(i). If $m\geq2$, then $d(u_2)+1=d(v_{m-1})+d(u_{n-1})$.
Obviously, $u_2=u_{n-1} $ when $n=3$ and $d(u_2)=d(u_{n-1})$ when
$n\geq4$ by Lemma \ref{p}~(i). Hence $d(v_{m-1})=1$, it contradicts
the fact that $d(v_{m-1})\geq d_{G_0}(v_{m-1})=2$. Therefore
$k\geq2$. Dually, we have $m\geq2$.

Note that $d(x_1)=d(x_3)=d(x_5)$, we have
$S(x_1)=2+d(w_1)=S(x_3)=2+d(v_1)=S(x_5)=2+d(u_1)$ by \eqref{s21} and
\eqref{sab}. It implies that $d(w_1)=d(v_1)=d(u_1)=3$ or $a+b$.

If $d(u_1)=3$. Applying \eqref{suv} with $(v,u)=(x_5,x)$, we have
$a=5-(a+b)$. So $a=2,b=1$. Furthermore, $S(u_1)=4+d(u_2)=7$ by
\eqref{s21} and \eqref{sab}. Thus $d(u_2)=3$. It follows from Lemma
\ref{p}~(i) that $d(u_i)=3$ for $2\leq i \leq n-1$. Similarly, we
have $d(v_i)=d(w_j)=3$ for  $2\leq i \leq m-1$ and $2\leq j \leq
k-1$. Note that $d_{G_0}(u_n)=a+b=3$. We have $d(u_n)=3$ and
$S(u_n)=7$ by \eqref{sab}. On the other hand, $S(u_n)=9$ by
\eqref{s21}, a contradiction.

If $d(u_1)=a+b\neq3$, then $a+b\geq4$ by Lemma \ref{3no}~(iii).
Applying \eqref{suv} with $(v,u)=(u_1,x)$, we have
$2=\frac{a+b-3+4+d(u_2)-(a+b)}{a+b-1}$. Note that $d_{G_0}(u_2)=2$
or $3$, we have $d(u_2)=2(a+b)-3>max\{a+b,d_{G_0}(u_2)\}$. It
contradicts Lemma \ref{3no}~(ii). $\square$
\begin{figure}[htbp]
\begin{center}
\includegraphics[totalheight=6.5cm]{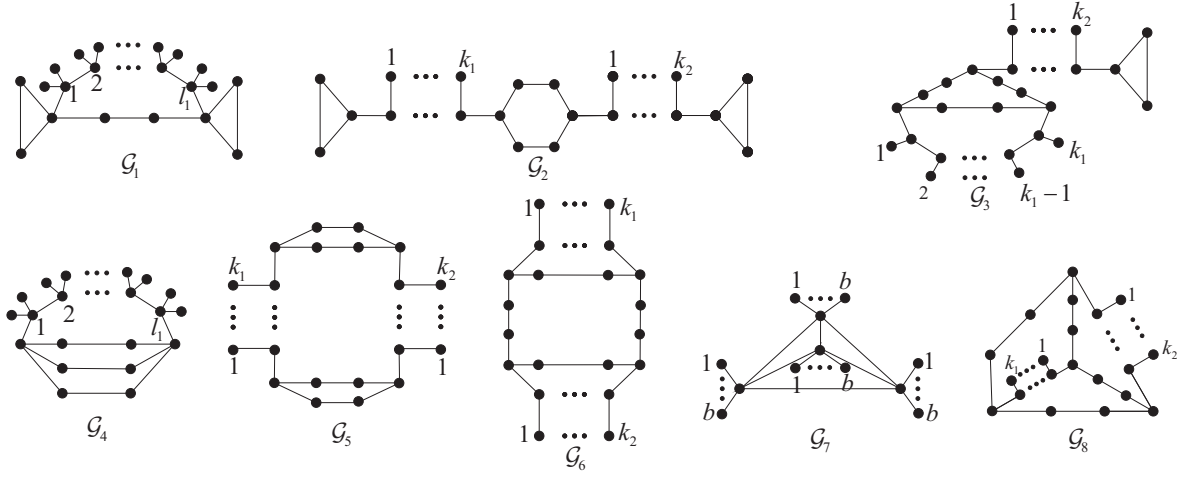}\\
\caption{\label{main1} \small{The graphs ${\mathcal G}_i$ for
i=1,\dots,8, where $l_1,b\geq1$ and $max\{k_1,k_2\}\geq1$ }}
\end{center}
\end{figure}
\begin{lemma}\label{t2} Let $G\in{\mathcal
G}$ with $G_0\in {\mathcal T}_2$ (see Fig.~1). Then $G=H_i$ for
$i=3,4,5,6$ (see Fig.~3) or $G\in{\mathcal G}_j$ for $j=1,2$ (see
Fig.~4).
\end{lemma}
\noindent {\bf Proof.}  If $G_0\in {\mathcal T}_2$, then
$d_{G_0}(w_1), d_{G_0}(r_1)\in \{3, 4\}$ and so each cycle of $G_0$
has the length of $3$ by Lemmas \ref{31} and \ref{p}~(iii). Hence
$G_0=T_2$ and $d(x_i)=2$ for $i=1,\dots,4$ (see Fig.~2), where
$k,l\geq1, n,m\geq 2$. For convenience, we set $w_k=v_1=u_1$ and
$v_m=u_n=r_l$. By \eqref{s21} and \eqref{sab},
$S(x_1)=2+d(w_1)=S(x_3)=2+d(r_1)$. Hence $d(w_1)=d(r_1)$.

If $G$ has no pendent vertex, then $k,l\in\{1,2,4\}$ and
$m,n\in\{2,4\}$ by Lemma \ref{31}.

First, let $k=1$. Then $d(w_1)=d(r_1)=4$. So $l=1$. Therefore
$G=H_i$ for $i=3,4$ (see Fig.~3). By \eqref{suv}, $H_3$ and $H_4$
are 2-walk $(2,2)$-linear and 2-walk $(1,4)$-linear, respectively.

Next, let $k=2$. Then $d(w_1)=d(r_1)=3$. By \eqref{s21} and
\eqref{sab}, $S(r_1)=4+d(r_2)= S(w_1)=7$. So $d(r_2)=3$ and $l=2$.
Similarly, $S(u_1)=3+d(u_2)+d(v_2)=S(w_1)=7$. So $d(u_2)=d(v_2)=2$.
Note that $n,m=2$ or $4$, we have $n=m=4$ and $d(u_3)=d(v_3)=2$.
Therefore, $G=H_5$ (see Fig.~3). By \eqref{suv}, $H_5$ is 2-walk
(2,1)-linear.

Finally, let $k=4$. With a similar argument of the case $k=2$, we
have $G=H_6$ is 2-walk (1,3)-linear (see Fig.~3).

If $G$ has at least one pendent vertex. We consider the following
two cases:

{\bf Case 1}. $k=l=1$. Then $d(w_1)=d(r_1)=4$ or $a+b$ by Lemma
\ref{3no}~(ii).

If $d(w_1)=4$. Applying \eqref{suv} with $(v,u)=(x_1,x)$, we have
$a=6-(a+b)$. It follows from Lemmas \ref{21} and \ref{3no}~(iii)
that $a+b=3$ or $4$.

If $a+b=3$, then $a=3,b=0$. So $S(w_1)=4+d(u_2)+d(v_2)=12$ by
\eqref{s21} and \eqref{sab}. By Lemma \ref{3no}~(ii),
$d(u_2),d(v_2)\in\{2,4,a+b\}$. Thus $d(u_2)=d(v_2)=4$. It implies
that $n=m=2$, which is impossible since $G$ is simple.

If $a+b=4$, then $a=b=2$. So $S(w_1)=4+d(u_2)+d(v_2)=10$ by
\eqref{s21} and \eqref{sab}. By Lemma \ref{3no}~(ii),
$d(u_2),d(v_2)\in\{2,4\}$. Without loss of generality, suppose that
$d(u_2)=2,d(v_2)=4$. Then $d(v_i)=4$ for $2\leq i \leq m-1$ by Lemma
\ref{p}~(i). For the vertex $u_2$, $S(u_2)=4+d(u_3)=6$ by
\eqref{s21} and \eqref{sab}. So $d(u_3)=2$. Thus $n\geq4$. Hence
$n=4$ by Lemma \ref{p}~(ii). Therefore $G\in {\mathcal G}_1$ (see
Fig.~4), where $l_1\geq1$. It is easy to see that any graph $G\in
{\mathcal G}_1$ is 2-walk $(2,2)$-linear.

If $d(w_1)=a+b>d_{G_0}(u_1)=4$. In this case, we show that there is
no such graph with exactly two main eigenvalues. Applying
\eqref{suv} with $(v,u)=(x_1,x)$ and $(v,u)=(w_1,x)$, respectively.
We have $a=2$ and $a=\frac{a+b-4+4+d(u_2)+d(v_2)-(a+b)}{a+b-1}$,
respectively. Hence $d(u_2)+d(v_2)=2(a+b)-2$. By Lemma
\ref{3no}~(ii) and the fact that $d(u_n)=d(r_1)=d(w_1)=a+b$, we have
$d(u_2), d(v_2) \in\{2,a+b\}$.

If $d(u_2)=2,d(v_2)\in\{2,a+b\}$ or $d(u_2)=a+b,d(u_2)=2$, then
$a+b=3$ or $4$. It contradicts the fact that $a+b>4$.

If $d(u_2)=d(v_2)=a+b$, then $2(a+b)-2=2(a+b)$, a contradiction.

{\bf Case 2}. $k\geq2$. Then $a=2$ by Lemma \ref{p}~(iii). By Lemmas
\ref{3no}~(ii), $d(w_1)=d(r_1)=3$ or $a+b$.

We first show that $d(w_1)=d(r_1)=3$. Otherwise, let
$d(w_1)=d(r_1)=a+b\neq3$. Then $a+b\geq4$ by Lemma \ref{3no}~(iii).
Applying \eqref{suv} with $(v,u)=(w_1,x)$, we have
$2=\frac{a+b-3+4+d(w_2)-(a+b)}{a+b-1}$. Note that $a+b\geq4$ and
$d_{G_0}(w_2)=2$ or $3$. We have
$d(w_2)=2(a+b)-3>max\{a+b,d_{G_0}(w_2)\}$. It contradicts Lemma
\ref{3no}~(ii). Hence $d(w_1)=d(r_1)=3$. It implies that $l\geq2$.

Next, applying \eqref{suv} with $(v,u)=(x_1,x)$, we have
$a=5-(a+b)$. So $a+b=3$ and $a=2,b=1$. For the vertex $w_1$,
$S(w_1)=4+d(w_2)=7$ by \eqref{s21} and \eqref{sab}. So $d(w_2)=3$.
Thus $d(w_i)=3$ for $2\leq i \leq k-1$ by Lemma \ref{p}~(i). Note
that $d_{G_0}(w_k)=3$. We have $d(w_k)=3$ by Lemma \ref{3no}~(ii).
Hence $d(w_i)=3$ for $1\leq i \leq k$.

Similarly, we have $d(r_j)=3$ for $1\leq j \leq l$.

For the vertex $w_k$, $S(w_k)=3+d(u_2)+d(v_2)=7$. Note that
$d(u_2),d(v_2)\in\{2,3\}$ by Lemma \ref{3no}~(ii). We have
$d(u_2)=d(v_2)=2$. It follows from Lemmas \ref{31} and \ref{p} that
$m=n=4$ and $d(u_3)=d(v_3)=2$.

Therefore $G\in {\mathcal G}_2$ (see Fig.~4), where
$max\{k_1,k_2\}\geq1$. It is easy to see that any graph $G\in
{\mathcal G}_2$ is 2-walk (2,1)-linear. $\square$

\begin{lemma}\label{t3} Let $G\in{\mathcal
G}$ with $G_0\in {\mathcal T}_3$ (see Fig.~1). Then $G=H_7$ (see
Fig.~3).
\end{lemma}
\noindent {\bf Proof.} If $G_0\in {\mathcal T}_3$, then
$d_{G_0}(r_1)\in \{3, 5\}$ and so the cycle of $G_0$ has the length
of $3$ by Lemmas \ref{31} and \ref{p}. Hence $G_0=T_3$ (see Fig.~2),
where $n,m,k\geq2,l\geq1$. For convenience, we set $u_1=v_1=w_1$ and
$u_n=v_m=w_k=r_l$.

We first show that $G$ contains at least one pendent vertex. On the
contrary, suppose that $G$ has no pendent vertex. Then
$n,m,k,l\leq4$ by Lemma \ref{31}. Without loss of generality,
suppose that $n\geq m\geq k$. If $l=1$, then $S(x_1)=7$ and
$S(u_2)=5$ or $8$ by \eqref{s21}. So $S(u_2)\neq S(x_1)$. On the
other hand, $S(u_2)=S(x_1)$ by \eqref{sab}, a contradiction. If
$l\geq2$, then $d(u_n)=4,d(r_1)=3$. So $S(x_1)=5$ and $S(u_{n-1})=6$
or $7$ by \eqref{s21}. On the other hand, $S(x_1)=S(u_{n-1})$ by
\eqref{sab}, a contradiction. Therefore, $G$ has at least one
pendent vertex. We consider the following two cases:

{\bf Case 1}. $l=1$. By Lemma \ref{p}~(iv), $d(x_{1})= d(x_{2})=3$
or $d(x_{1})= d(x_{2})=2$.

If $d(x_{1})= d(x_{2})=3$, then $a=3,b=0,d(u_n)=5$ by Lemma
\ref{p}~(iv). So $d(u_1)=3\neq d(u_n)$ by Lemma \ref{3no}~(ii). It
follows from Lemma \ref{p}~(i) that $n,m,k\leq3$. Without loss of
generality, suppose that $n=m=3,k=2$ or $3$. Then
$S(u_3)=6+d(u_2)+d(v_2)+d(w_{k-1})=15$ by \eqref{s21} and
\eqref{sab}.  Note that $d(u_2),d(v_2),d(w_{k-1})=2$ or $3$ by Lemma
\ref{3no}~(ii). We have $d(u_2)=d(v_2)=d(w_{k-1})=3$. So
$S(u_1)=6+d(w_2)=9$ by \eqref{s21} and \eqref{sab}. Thus $d(w_2)=3$.
It implies that $k=3$. Therefore $G=H_7$ (see Fig.~3). By
\eqref{suv}, $H_7$ is 2-walk $(3,0)$-linear.

If $d(x_{1})= d(x_{2})=2$. We show that in this case there is no
such graph with exactly two main eigenvalues. We consider the
following two cases:

{\bf Subcase 1}. $max\{n,m,k\}\geq4$. Without loss of generality,
suppose that $n\geq 4$. Then $d(u_1)=d(u_{n})$ and $d(u_i)=d(u_{2})$
for $2\leq i\leq n-1$ by Lemma \ref{p}~(i). Note that
$d_{G_0}(u_1)\neq d_{G_0}(u_{n})$. We have $d(u_1)=d(u_{n})=a+b\geq
d_{G_0}(u_{n})=5$ by Lemma \ref{3no}~(ii). Hence
\begin{equation*}
S(u_1)=a+b-3+d(u_2) +d(v_2)+d(w_2),\
S(u_n)=a+b-5+4+d(u_{n-1})+d(v_{m-1})+d(w_{k-1}).
\end{equation*}

We next show that $d(v_2)=d(v_{m-1})$. If $m=2$, then $v_2=u_n$,
$v_{m-1}=u_1$. So $d(v_2)=d(u_n)=d(u_1)=d(v_{m-1})$. If $m=3$,
clearly, $d(v_2)=d(v_{m-1})$. If $m\geq 4$, then $d(v_2)=d(v_{m-1})$
by Lemma \ref{p}~(i). Therefore $d(v_2)=d(v_{m-1})$ for all
$m\geq2$.

Similarly, we have $d(w_2)=d(w_{k-1})$ for all $k\geq2$. Hence
$S(u_1)\neq S(u_n)$. On the other hand, $S(u_1)=S(u_n)$ by
\eqref{sab}, a contradiction.

{\bf Subcase 2}. $max\{n,m,k\}\leq3$. Without loss of generality,
suppose that $n=m=3$ and $k=2$ or $3$. We claim that $d(u_2)=a+b$.
Otherwise, let $d(u_2)=d_{G_0}(u_2)=2$. Then
$d(u_1)+d(u_3)=S(u_2)=S(x_1)=2+d(u_3)$ by \eqref{s21} and
\eqref{sab}, which is impossible since $d(u_1)\geq d_{G_0}(u_1)=3$.
Hence $d(u_2)=a+b$. Similarly, we have $d(v_2)=a+b$.

Note that $d(u_3)\in\{a+b,5\}$ by Lemma \ref{3no}~(ii). We consider
the following two cases:

If $d(u_3)=a+b\geq d_{G_0}(u_3)=5$.  Applying Lemma \ref{3no}~(iv)
with $C=u_1u_2u_3v_2u_1$, we have $d(u_1)\neq a+b$. So $d(u_1)=3$.
 Applying \eqref{suv} with $(v,u)=(u_1,x)$ and
$(v,u)=(x_1,x)$, respectively. We have $a=\frac{a+b+d(w_2)}{2}$ and
$a=2$, respectively. Hence $d(w_2)=4-(a+b)<0$, a contradiction.

If $d(u_3)=d_{G_0}(u_3)=5\neq a+b$, then $a=7-(a+b)$ by applying
\eqref{suv} with $(v,u)=(x_1,x)$. It follows from Lemma
\ref{3no}~(iii) that $a+b=3$ or $4$. If $a+b=3$, then $a=4,b=-1$.
Hence $S(u_2)=d(u_1)+6=11$ by \eqref{s21} and \eqref{sab}. This is
impossible since $d(u_1)=3$ by Lemma \ref{3no}~(ii). If $a+b=4$,
then $a=3,b=1$. Thus $S(u_2)=d(u_1)+7=13$ by \eqref{s21} and
\eqref{sab}. This is also impossible since $d(u_1)=3$ or $4$.

{\bf Case 2}. $l\geq2$. We show that in this case there is no such
graph with exactly two main eigenvalues. By Lemma \ref{p}~(iii), we
have $a=2$ and $d(x_{1})=d(x_2)= 2$. Applying \eqref{suv} with
$(v,u)=(x_1,x)$, we have $2=2+d(r_1)-(a+b)$. So $d(r_1)=a+b$.
Applying \eqref{suv} with $(v,u)=(r_1,x)$, we have
$2=\frac{a+b-3+4+d(r_2)-(a+b)}{a+b-1}$. Hence $d(r_2)=2(a+b)-3$. By
Lemma \ref{3no}~(ii), $d(r_2)\in\{2,4,a+b\}$.

If $d(r_2)=2$ or $4$, then $a+b$ is not an integer, a contradiction.

If $d(r_2)=a+b$, then $a+b=3$. So $a=2,b=1$. By Lemmas
\ref{3no}~(ii) and \ref{p}~(i), we have $d(r_l)=4$ and
$d(r_{l-1})=d(r_2)=3$. Hence
$S(r_{l})=3+d(u_{n-1})+d(w_{k-1})+d(v_{m-1})=9$ by \eqref{s21} and
\eqref{sab}. It follows that $d(u_{n-1})=d(w_{k-1})=d(v_{m-1})=2$.
Thus $S(u_{n-1})=5$ by \eqref{sab}. On the other hand,
$S(u_{n-1})=4+d(u_{n-2})\geq6$ by \eqref{s21}, a contradiction.
$\square$

\begin{lemma}\label{t4} Let $G\in{\mathcal
G}$ with $G_0\in {\mathcal T}_4$ (see Fig.~1). Then $G=H_i$ for
$i=8,9,10$ (see Fig.~3) or $G\in {\mathcal G}_3$ (see Fig.~4).
\end{lemma}
\noindent {\bf Proof.} If $G_0\in {\mathcal T}_4$, then
$d_{G_0}(w_1)\in \{3, 4\}$ and so the cycle of $G_0$ has the length
of $3$ by Lemmas \ref{31} and \ref{p}~(iii). Hence $G_0=T_4$ (see
Fig.~2), where $k\geq1$ and $ p,q,n,m\geq2$. For convenience, we set
$u_1=v_1=s_1,u_n=v_m=t_1$ and $s_p=t_q=w_k$.

If $G$ has no pendent vertex. Without loss of generality, we may
assume that $n\geq m,p\geq q$ and consider the following two cases:

{\bf Case 1}. $k=1$. Then $S(x_1)=6$  by \eqref{s21}. We claim that
$p=2$. Otherwise, let $p\geq 3$. Then $d(s_2)=2$ and $S(s_2)=5$ or
$7$ by \eqref{s21}. So $S(x_1)\neq S(s_2)$. On the other hand,
$S(x_1)=S(s_2)$  by \eqref{sab}, a contradiction. Hence $p=2$.
Similarly, we have $q=2,n=3,m=2$ or $3$. If $m=2$, then applying
\eqref{suv} with $(v,u)=(w_1,u_2)$ and $(v,u)=(w_1,u_1)$,
respectively, we have $a=2$ and $a=1$, respectively. A
contradiction. If $m=3$, then $G=H_8$ (see Fig.~3). By \eqref{suv},
$H_8$ is 2-walk $(2,2)$-linear.

{\bf Case 2}. $k>1$. By Lemma \ref{31}, we have $k,p,q,m,n=2$ or
$4$.

If $k=2$, then $S(w_1)=7=S(w_2)=3+d(s_{p-1})+d(t_{q-1})$ by
\eqref{s21} and \eqref{sab}. So $d(s_2)=d(t_2)=2$. It implies that
$p=q=4$. Hence $n=4,m=2$. Otherwise, let $m=n=4$. Then $S(u_1)=6\neq
S(w_1)=7$ by \eqref{s21}. On the other hand, $S(u_1)=S(w_1)$ by
\eqref{sab}, a contradiction. Therefore $n=4,m=2$ and $G=H_9$ (see
Fig.~3). By \eqref{suv}, $H_9$ is 2-walk $(2,1)$-linear.

If $k=4$, with a similar argument, we have $G=H_{10}$ is 2-walk
$(1,3)$-linear (see Fig.~3).

If $G$ has at least one pendent vertex, then $k>1$. Otherwise,
suppose that $k=1$. Then $d(w_1)=4$ or $a+b$ by Lemma
\ref{3no}~(ii).

If $d(w_1)=4\neq a+b$. Applying \eqref{suv} with $(v,u)=(x_1,x)$, we
have $a=6-(a+b)$. It follows from Lemma \ref{3no}~(iii) that
$a=3,b=0$. So $S(w_1)=4+ d(s_{p-1})+d(t_{q-1})=12$ by  \eqref{s21}
and \eqref{sab}. This is impossible since $d(s_{p-1}),d(t_{q-1})=2$
or $3$ by Lemma \ref{3no}~(ii).

If $d(w_1)=a+b\geq 4$. Applying \eqref{suv} with $(v,u)=(x_1,x)$ and
$(v,u)=(w_1,x)$, respectively, we have $a=2$ and
$a=\frac{a+b-4+4+d(t_{q-1})+d(s_{p-1})-(a+b)}{a+b-1}$, respectively.
Thus
\begin{equation}\label{9a}
d(t_{q-1})+d(s_{p-1})=2(a+b)-2.
\end{equation}
Note that $d(t_{q-1}),d(s_{p-1})\in\{2,3,a+b\}$. We consider the
following five cases by symmetry:

If $d(t_{q-1})=d(s_{p-1})=2$, then $a+b=3$. This contradicts the
fact that $a+b\geq 4$.

If $d(t_{q-1})=2,d(s_{p-1})=3$, then $a+b=3.5$. This contradicts
Lemma \ref{21}.

If $d(t_{q-1})=2,d(s_{p-1})=a+b $, then $a+b=4$ by \eqref{9a}. Note
that $a=2$. We have $b=2$ and $d(s_{p-1})=d(w_1)=4$. By Lemma
\ref{p}, $d(s_{i})=4$ for $2\leq i \leq p-1$. So
$S(s_2)=6+d(u_1)=10$ by \eqref{s21} and \eqref{sab}. Thus
$d(u_1)=4$. Similarly, we have $d(u_n)=4$. By \eqref{s21} and
\eqref{sab}, $S(u_1)=5+d(u_2)+d(v_2)=10$. This is impossible since
$d(u_2),d(v_2)=2$ or $4$.

If $d(t_{q-1})=3,d(s_{p-1})=a+b$, then $a+b=5$ by \eqref{9a}. Recall
that $a=2$, we have $b=3$ and $d(w_1)=5$. Note that
$d(t_{q-1})=3\neq a+b$. We have $q=2$. Thus
$S(t_{q-1})=9=5+d(u_{n-1})+d(v_{n-1})$ by \eqref{s21} and
\eqref{sab}. By Lemma \ref{3no}~(ii),
$d(u_{n-1}),d(v_{n-1})\in\{2,3,5\}$. Hence
$d(u_{n-1})=d(v_{n-1})=2$. Therefore $S(u_{n-1})=3+d(u_{n-2})=7$,
which is impossible since $d(u_{n-2})\in\{2,3,5\}$.

If $d(t_{q-1})=d(s_{p-1})=a+b$, then $2(a+b)=2(a+b)-2$ by
\eqref{9a}, a contradiction.

Hence $k>1$. It follows from Lemma \ref{p}~(iii) that $a=2$ and $
d(x_1)=d(x_2)=2$.

Applying \eqref{suv} with $(v,u)=(x_1,x)$, we have
$2=2+d(w_1)-(a+b)$. So $d(w_1)=a+b$. Applying \eqref{suv} with
$(v,u)=(w_1,x)$, we have $2=\frac{a+b-3+4+d(w_2)-(a+b)}{a+b-1}$.
Thus $d(w_2)=2(a+b)-3\geq3$. Note that $d(w_2)\in\{2,3,a+b\}$ by
Lemma \ref{3no}~(ii). We have $d(w_2)=3$ or $a+b$. It follows that
$a+b=3$ and $a=2,b=1$.

By Lemma \ref{p}~(i), $d(w_i)=d(w_2)=3$ for $2\leq i \leq k-1$.

By Lemma \ref{3no}~(ii), $d(w_1)=d(w_k)=d(u_1)=d(u_n)=3$.

For the vertex $w_k$, we have $S(w_k)=3+d(s_{p-1})+d(t_{q-1})=7$. It
follows from Lemma \ref{3no}~(ii) that $d(s_{p-1})=d(t_{q-1})=2$ and
$p,q\geq3$. Hence $p=q=4$ and $d(s_2)=d(t_2)=2$ by Lemmas \ref{31}
and \ref{p}~(ii).

For the vertex $u_1$, we have $S(u_1)=2+d(v_2)+d(u_2)=7$. Note that
$d(u_2),d(v_2)\in\{2,3\}$. We may assume that $d(v_2)=2,d(u_2)=3$ by
symmetry. It follows from Lemmas \ref{31} and \ref{p}~(ii) that
$m=4,d(v_3)=2$ and $d(u_i)=3$ for $2\leq i\leq n-1$.

Therefore $G\in{\mathcal G}_3$ (see Fig.~4), where
$max\{k_1,k_2\}\geq1$. It is easy to see that any graph
$G\in{\mathcal G}_3$ is 2-walk $(2,1)$-linear.

Up to now, we complete the proof of the Lemma. $\square$

\begin{lemma}\label{t5} There is no graph $G\in{\mathcal
G}$ with $G_0\in{\mathcal T}_{5}$ (see Fig.~1).
\end{lemma}
\noindent {\bf Proof.}  By way of contradiction, suppose that
$G\in{\mathcal G}$  with $G_0\in{\mathcal T}_{5}$. Let $G_0=T_{5}$
(see Fig.~2), where $n,m,k,p,q\geq 2$. For convenience, we set
$u_1=v_1=s_1,u_n=w_k=t_q, v_m=s_p=w_1=t_1$ and consider the
following two cases:

{\bf Case 1}.  There is a vertex $v\in\{v_{m-1},s_{p-1},w_2,t_2\}$
such that $d_{G_0}(v)=2$ and $d(v)=a+b$. Without loss of generality,
let $v=v_{m-1}$. Then $d(v_i)=a+b$ for $2\leq i \leq m-1$ by Lemma
\ref{p}~(i). In particular, $a\geq 2,a+b\geq3$ by Lemma
\ref{3no}~(iii). We first show that $d(u_1)=d(w_1)=a+b$.

If  $m\geq 4$, then $d(u_1)=d(w_1)$ by Lemma \ref{p}~(i). Note that
$d_{G_0}(u_1)\neq d_{G_0}(w_1)$. We have $d(u_1)=d(w_1)=a+b$ by
Lemma \ref{3no}~(ii).

If $m=3$. Applying \eqref{suv} with $(v,u)=(v_2,x)$, we have
\begin{equation*}
a=\frac{a+b-2+d(u_1)+d(w_1)-(a+b)}{a+b-1}=\frac{d(u_1)+d(w_1)-2}{a+b-1}.
\end{equation*}
By Lemma \ref{3no}~(ii), $d(u_1)=3$ or $a+b$ and $d(w_1)=4$ or
$a+b$.

If $d(u_1)=3, d(w_1)=4$, then $a=\frac{5}{a+b-1}$. Note that $a\geq
2$ is an integer. We have $a+b=2$. It contradicts the fact that
$a+b\geq3$.

If $d(u_1)=3,d(w_1)=a+b\geq 4$, then $a=1+\frac{2}{a+b-1}<2$, a
contradiction.

If $d(u_1)=a+b\geq3,d(w_1)=4\neq a+b$, then $a=1+\frac{3}{a+b-1}$ is
not an integer, a contradiction.

Hence $d(u_1)=d(w_1)= a+b$.

Applying Lemma \ref{3no}~(iv) with $C=u_1v_2\dots v_{m-1}
w_1s_{p-1}\dots s_2u_1$, we have $p\geq3$ and $d(s_i)\neq a+b$ for
some $2\leq i \leq p-1$. It follows from Lemmas \ref{3no}~(ii) and
\ref{p}~(i) that $d(s_i)=2$ for $2\leq i \leq p-1$. Applying
\eqref{suv} with $(v,u)=(v_2,x)$ and $(v,u)=(u_1,x)$, respectively.
We have $a=2$ and $a=1+\frac{d(u_2)}{a+b-1}$, respectively. This
together with Lemma \ref{3no}~(ii) and the fact that $a+b\geq
d_{G_0}(w_1)=4$ implies that $d(u_2)=a+b-1=d_{G_0}(u_2)\geq3$. Hence
$n=2, d(u_2)=d_{G_0}(u_2)=3$. Therefore $a+b=4$ and $a=b=2$. By
\eqref{s21} and \eqref{sab}, $S(u_2)=4+d(w_{k-1})+d(t_{q-1})=8$. By
Lemma \ref{3no}~(ii), $d(w_{k-1}),d(t_{q-1})=2$ or $4$. Thus
$d(w_{k-1})=d(t_{q-1})=2$. Hence $S(w_{k-1})=3+d(w_{k-2})=6$ by
\eqref{s21} and \eqref{sab}, which is impossible since
$d(w_{k-2})=2$ or $4$.

{\bf Case 2}.  For any vertex $v\in\{v_{m-1},s_{p-1},w_2,t_2\}$, we
have $d_{G_0}(v)=3$ or  $d(v)=2$. It follows from Lemma \ref{31}
that $m,k,p,q\leq4$. Without loss of generality, suppose that $m\geq
p,k\geq q$. Hence $m,k=3$ or $4$.

{\bf Subcase 1}. $max\{m,k\}=4$. Without loss of generality, suppose
that $m=4$. Then $d(u_1)=d(w_1)$ by Lemma \ref{p}~(i). Note that
$d_{G_0}(u_1)\neq d_{G_0}(w_1)$. We have $d(u_1)=d(w_1)=a+b\geq
d_{G_0}(w_1)=4$ by Lemma \ref{3no}~(ii). Thus $G$ has at least one
pendent vertex. Applying \eqref{suv} with $(v,u)=(v_2,x)$ and
$(v,u)=(u_1,x)$, respectively. We have $a=2$ and
$a=\frac{a+b-3+2+d(s_2)+d(u_2)-(a+b)}{a+b-1}$, respectively. Hence
$d(s_2)+d(u_2)=2(a+b)-1$.

If $p>2$, then $d(s_2)=2$. It follows from the fact that
$d_{G_0}(u_2)=2$ or $3$ and $a+b\geq 4$ that $d(u_2)=2(a+b)-3>
max\{a+b,d_{G_0}(u_2)\}$, which is impossible by Lemma
\ref{3no}~(ii).

If $p=2$, then $d(s_2)=d(w_1)=a+b$. So $d(u_2)=a+b-1\geq 3$. It
follows that $n=2,d(u_2)=d_{G_0}(u_2)=3$ and $a+b=4$. Hence $a=b=2$.
By \eqref{sab}, $S(w_{k-1})=6$. On the other hand,
$S(w_{k-1})=d(w_{k-2})+d(u_2)=5$ or $7$ by \eqref{s21}, a
contradiction.

{\bf Subcase 2}. $m=k=3$. Then
$d(w_1)+d(u_1)=S(v_2)=S(w_2)=d(w_1)+d(u_n)$ by \eqref{s21} and
\eqref{sab}. So $d(u_1)=d(u_n)$. Thus $S(u_1)=S(u_n)$ by
\eqref{sab}.

We claim that $p=q=2$ or $3$. Otherwise, let $p\neq q$. Without loss
of generality, suppose that $p=3,q=2$. Then
$S(u_1)=d(u_1)-3+4+d(u_2)$ and
$S(u_n)=d(u_n)-3+2+d(w_1)+d(u_{n-1})$. Note  that
$d(u_2)=d(u_{n-1}),d(w_1)>2$. We have $S(u_1)\neq S(u_n)$, a
contradiction. Hence $p=q=2$ or $3$. We consider the following two
cases:

{\bf Subcase 2.1}. $G$ has no pendent vertex. Then $n=2$. Otherwise,
suppose that $n>2$. Then $d(u_2)=d(v_2)=2$. By \eqref{s21} and
\eqref{sab}, $3+d(u_3)=S(u_2)=S(v_2)=7$. This is impossible since
$d(u_3)=2$ or $3$. Hence $n=2$. If $p=q=2$. Applying \eqref{suv}
with $(v,u)=(u_1,v_2)$ and $(v,u)=(w_1,u_1)$, respectively. We have
$a=2$ and $a=1$, respectively. A contradiction. If $p=q=3$, also
applying \eqref{suv} with $(v,u)=(u_1,v_2)$ and $(v,u)=(w_1,u_1)$,
respectively, we have $a=0$ and $a=1$, respectively. Also a
contradiction.

{\bf Subcase 2.2}. $G$ has at least one pendent vertex. Then
$a\geq2,a+b\geq3$.

First, let $p=q=3$. Applying \eqref{suv} with $(v,u)=(w_1,x)$, we
have $a=\frac{d(w_1)-4+8-(a+b)}{d(w_1)-1}$. By Lemma \ref{3no}~(ii),
$d(w_1)=a+b$ or $4$. It follows that $a<2$, a contradiction.

Next, suppose that $p=q=2$. By Lemma \ref{3no}~(ii), $d(w_1)=a+b$ or
$4$.

If $d(w_1)=a+b\geq4$. Applying \eqref{suv} with $(v,u)=(w_1,x)$ and
note that $d(u_1)=d(u_n)$, we have $a=\frac{2d(u_1)}{a+b-1}$. If
$d(u_1)=a+b$, then $a=2+\frac{2}{a+b-1}$ is not an integer, which is
impossible by Lemma \ref{21}. If $d(u_1)=d_{G_0}(u_1)=3\neq a+b$,
then $a+b=4$ and $a=b=2$. So $S(v_2)=6$ by \eqref{sab}. On the other
hand, $S(v_2)=d(u_1)+d(w_1)=7$ by \eqref{s21}, a contradiction.

If $d(w_1)=4\neq a+b$. Applying \eqref{suv} with $(v,u)=(v_2,x)$, we
have $a=4+d(u_1)-(a+b)$. If $d(u_1)=a+b$, then $a=4$. For the vertex
$u_1$, we have $S(u_1)=4(a+b)+b=a+b-3+6+d(u_2)$. It follows that
$d(u_2)=4(a+b)-7>max\{a+b,d_{G_0}(u_2)\}$, which is impossible by
Lemma \ref{3no}~(ii). If $d(u_1)=3\neq a+b$, then $a=7-(a+b)$. Note
that $a+b\neq3,4$. We have $a+b=5$ and $a=2,b=3$ by Lemma
\ref{3no}~(iii). By \eqref{s21} and \eqref{sab},
$S(w_1)=11=4+d(u_1)+d(u_n)$. This is impossible since
$d(u_1)=d(u_n)$.

Up to now, we have completed the proof of the Lemma. $\square$.

\begin{lemma}\label{t6} Let $G\in{\mathcal
G}$ with $G_0\in {\mathcal T}_{6}$. Then $G=H_i$ for $i=11,\dots,15$
(see Fig.~3) or $G\in{\mathcal G}_4$ (see Fig.~4).
\end{lemma}
\noindent {\bf Proof.} Let $G_0=T_{6}$ (see Fig.~2), where
$n,m,p,k\geq2$. For convenience, we set $u_1=v_1=w_1=s_1$ and
$u_n=v_m=w_k=s_p$.

{\bf Case 1}. There is a vertex $v\in\{u_2,v_{2},w_2,s_{2}\}$ such
that $d_{G_0}(v)=2$ and $d(v)=a+b$. Without loss of generality,
suppose that $v=u_{2}$. Applying \eqref{suv} with $(v,u)=(u_2,x)$,
we have
\begin{equation}\label{666}
a=\frac{a+b-2+d(u_1)+d(u_3)-(a+b)}{a+b-1}=\frac{d(u_1)+d(u_3)-2}{a+b-1}.
\end{equation}
By Lemmas \ref{3no}~(ii), $d(u_1),d(u_3)=4$ or $a+b$. We consider
the following three cases:

{\bf Subcase 1}. $d(u_1)=d(u_3)=a+b$, then $a=2$ by \eqref{666}.

We now show that $d(u_i)=a+b$ for $1\leq i\leq n$. If $n=3$, then
obviously, $d(u_i)=a+b$ for $1\leq i\leq n$. If $n\geq4$, then
$d(u_i)=d(u_2)=a+b$  for $1< i<n$ and $d(u_n)=d(u_1)=a+b$ by Lemma
\ref{p}~(i). Hence $d(u_i)=a+b$ for $1\leq i\leq n$.

Applying Lemma \ref{3no}~(iv) with $C=u_1u_2\dots u_nv_{m-1}\dots
v_{2}u_1$, we have $d(v_i)\neq a+b$ for some $2\leq i\leq m-1$. It
follows from Lemmas \ref{3no}~(ii) and \ref{p}~(i) that $m\geq3$ and
$d(v_i)=2$ for all $2\leq i\leq m-1$. Thus $S(v_2)=2a+b=a+b+d(v_3)$.
Note that $a=2$. We have $d(v_3)=2$. So $m\geq 4$. It follows from
Lemma \ref{p}~(ii) that $m=4$. Similarly, we have $k=p=4$ and
$d(w_2)=d(w_3)=d(s_2)=d(s_3)=2$.

For the vertex $u_1$, $S(u_1)=2(a+b)+b=a+b-4+6+a+b$ by \eqref{s21}
and \eqref{sab}. Hence $b=2$. It follows that $d(u_i)=a+b=4$ for
$1\leq i \leq n$.

Therefore $G\in {\mathcal G}_4$ (see Fig.~4), where $l_1\geq1$. It
is easy to see that any graph $G\in {\mathcal G}_4$ is 2-walk
$(2,2)$-linear.

{\bf Subcase 2}. $d(u_1)=a+b,d(u_3)=4$ or $d(u_1)=4,d(u_3)=a+b$ and
$a+b\neq4$. Then $a=1+\frac{3}{a+b-1}$ is not an integer by
\eqref{666}, a contradiction.

{\bf Subcase 3}. $d(u_1)=d(u_3)=4\neq a+b$. Then $n=3$ and
$a=\frac{6}{a+b-1}$. It follows from Lemma \ref{3no}~(iii) that
$a+b=3$ and $a=3,b=0$. Thus $d(u_2)=a+b=3$. For the vertex $u_1$,
$S(u_1)=12=3+d(v_2)+d(w_2)+d(s_2)$ by \eqref{s21} and \eqref{sab}.
Note that $d(v_2),d(w_2),d(s_2)\in\{2,4,a+b\}$. We have
$d(v_2)=d(w_2)=d(s_2)=a+b=3$. So $S(v_2)=9=5+d(v_3)$. Thus
$d(v_3)=4\neq a+b$. It implies that $m=3$. Similarly, we have
$k=p=3$. Therefore $G=H_{11}$ (see Fig.~3). By \eqref{suv}, $H_{11}
$ is 2-walk $(3,0)$-linear.

{\bf Case 2}.  For any vertex  $v\in\{u_2,v_{2},w_2,s_{2}\}$,
$d_{G_0}(v)=4$ or  $d(v)=2$. By Lemma \ref{31}, $n,m,p,k\leq4$.
Without loss of generality, suppose that $n\geq m\geq p\geq k$. Then
$n\geq m\geq p\geq3$ and $d(u_2)=d(v_2)=d(s_2)=2$. By \eqref{s21}
and \eqref{sab}, $S(u_2)=d(u_1)+d(u_3)=S(v_2)=d(u_1)+d(v_3)$. Thus
$d(v_3)=d(u_3)$. It implies that $m=n$. Similarly, we have $p=n$ and
$k=n$ or $2$. Hence $k=p=m=n\in\{3,4\}$ or $k=2,p=m=n\in\{3,4\}$.

If $G$ has no pendent vertex, then  $G=H_{i}$ for $i=12,\dots,15$
(see Fig.~3). By \eqref{suv},  $H_{12}$ is 2-walk $(1,6)$-linear,
$H_{13}$ is 2-walk $(0,8)$-linear, $H_{14}$ is 2-walk $(2,2)$-linear
and $H_{15}$ is 2-walk $(1,4)$-linear.

If $G$ has at least one pendent vertex $x$, then $x\in N(u_1)$ or
$x\in N(u_n)$. Without loss of generality, suppose that $x\in
N(u_1)$. Then $d(u_1)=a+b\geq5$.  Applying \eqref{suv} with
$(v,u)=(u_1,x)$, we have $a=\frac{d(w_2)+2}{a+b-1}$. By Lemma
\ref{3no}~(ii), $d(w_2)=2,4$ or $a+b$. It implies that $a<2$. This
is impossible by Lemma \ref{3no}~(iii).

Up to now, we have competed the proof of the Lemma. $\square$

\begin{lemma}\label{t7} Let $G\in{\mathcal
G}$ with  $G_0\in {\mathcal T}_{7}$(see Fig.~1). Then $G=H_{i}$ for
$i=16,\dots,23$ (see Fig.~3) or $G\in{\mathcal G}_j$ for $j=5,6$
(see Fig.~4).
\end{lemma}
\noindent {\bf Proof.} Let $G_0=T_{7}$ (see Fig.~2), where
$n,m,k,l,p,q\geq2$. For convenience, we set $u_1=w_1=s_1,
u_n=w_k=t_1, v_1=r_1=s_{p}$ and $v_m=r_l=t_q$.

If $G$ has no pendent vertex, then $n,m,k,l,p,q\leq4$ by Lemma
\ref{31}. Without loss of generality, suppose that $n\geq k, m\geq
l, p\geq q$. Then $n=3$ or $4$. By \eqref{s21} and \eqref{sab},
$S(u_1)=d(u_2)+d(w_2)+d(s_2)=S(u_n)=d(u_{n-1})+d(w_{k-1})+d(t_2)$.
Note that $d(u_2)=d(u_{n-1}), d(w_2)=d(w_{k-1})$. We have
$d(s_2)=d(t_2)$. It implies that $p=q$. Similarly, we have $k=l$. If
$n=3$, then $m=3,k,p=2$ or $3$ by Lemma \ref{31}. Hence $G=H_i$ for
$i=16,17,18,19$ (see Fig.~3). If $n=4$, then $m=4,k,p=2$ or $4$ by
Lemma \ref{31}. Hence $G=H_j$ for $j=20,21,22,23$ (see Fig.~3). By
\eqref{suv}, $H_{16}$ is 2-walk $(2,2)$-linear, $H_{17}$ and
$H_{18}$ are 2-walk $(1,4)$-linear, $H_{19}$ is 2-walk
$(0,6)$-linear, $H_{20}$ is 2-walk $(3,-1)$-linear, $H_{21}$ and
$H_{22}$ are 2-walk $(2,1)$-linear, $H_{23}$ is 2-walk
$(1,3)$-linear.

If $G$ has at least one pendent vertex. Then $a\geq2,a+b\geq3$.

We first show that $d(v_1)=d(v_m)$.

If $max\{m,l\}\geq4$, then $d(v_1)=d(v_m)$ by Lemma \ref{p}~(i).

If $m,l\leq3$. By way of contradiction, suppose that $d(v_1)\neq
d(v_m)$. By Lemma \ref{3no}~(ii), we may assume that $d(v_1)=a+b>3,
d(v_m)=3$ without loss of generality. Applying \eqref{suv} with
$(v,u)=(v_1,v_3)$, we have
\begin{equation*}
a=
\begin{cases}
\frac{d(s_{p-1})-d(t_{q-1})}{a+b-3} &\mbox{if}\  m=3,l=2\ \mbox{or}\ m=2,l=3,\\
1+\frac{d(s_{p-1})-d(t_{q-1})}{a+b-3}\quad& \mbox{if}\  m=3,l=3.
\end{cases}
\end{equation*}
By Lemma \ref{3no}~(ii), $d(s_{p-1}),d(t_{q-1})=2,3$ or $a+b$.

If $m=3,l=2$ or $m=2,l=3$, then $d(s_{p-1})=a+b,d(t_{q-1})=2$ since
$a\geq2,a+b\geq3$. So $a=1+\frac{1}{a+b-3}$. It implies that $a+b=4$
and $a=b=2$. By \eqref{s21} and \eqref{sab}, $S(v_3)=6+d(v_2)=8$. So
$d(v_2)=2$ and $S(v_2)=6$ by \eqref{sab}. On the other hand,
$S(v_2)=7$ by \eqref{s21}, a contradiction.

If $m=l=3$, then $d(s_{p-1})=a+b\neq 3,d(t_{q-1})=3$ or
$d(s_{p-1})=a+b,d(t_{q-1})=2$ since $a\geq2,a+b\geq3$. If
$d(s_{p-1})=a+b\neq 3, d(t_{q-1})=3$, then $a+b\geq4$ and $a=2$.
Thus $b\geq2$. By \eqref{s21} and \eqref{sab},
$S(v_3)=6+b=d(v_2)+d(r_2)+3$, which is impossible since
$d(v_2),d(r_2)=2$ or $2+b$ and $b\geq2$. If
$d(s_{p-1})=a+b,d(t_{q-1})=2$, then $a=2+\frac{1}{a+b-3}$. It
implies that $a+b=4$ and $a=3,b=1$. By \eqref{s21} and \eqref{sab},
$S(v_3)=d(v_2)+d(r_2)+2=10$. Note that $d(v_2),d(r_2)=2$ or $4$ by
Lemma \ref{3no}~(ii), we have $d(v_2)=d(r_2)=4$.  Thus $S(v_{2})=13$
by \eqref{sab}. On the other hand, $S(v_{2})=9$ by \eqref{s21}, a
contradiction.

Hence $d(v_1)=d(v_m)$.

Next, we show that $d(v_1)=d(v_m)=a+b$. On the contrary, suppose
that $d(v_1)=d(v_m)=3\neq a+b$ by Lemma \ref{3no}~(ii). Then
$a+b\geq4 $. Note that $max\{m,l\}\geq3$. Without loss of
generality, suppose that $m\geq3$. By Lemma \ref{3no}~(ii), we have
$d(v_2)=a+b$ or $2$.

If $d(v_2)=a+b$. Applying \eqref{suv} with $(v,u)=(v_2,x)$, we have
$a=\frac{d(v_3)+1}{a+b-1}$. By Lemma \ref{3no}~(ii), $d(v_3)=2,3$ or
$a+b$. This together with $a+b\geq4$ implies that $a<2$. It
contradicts the fact that $a\geq2$.

If $d(v_2)=2$. Applying \eqref{suv} with $(v,u)=(v_2,x)$, we have
$a=3+d(v_3)-(a+b)$. If $m\geq4$, then $d(v_3)=2$ by Lemma
\ref{p}~(i). Note that $a+b\geq4$. We have $a=5-(a+b)\leq1$. It
contradicts the fact $a\geq2$. Thus $m=3$. Hence $d(v_3)=3$ and
$a=6-(a+b)$. Note that $a+b\geq4,a\geq2$. We have $a+b=4$ and
$a=b=2$. By \eqref{s21} and \eqref{sab},
$S(v_1)=2+d(r_2)+d(s_{p-1})=8$. So $d(r_2)+d(s_{p-1})=6$.

If $l=2$, then $d(r_2)=d(v_m)=3$. So $d(s_{p-1})=3\neq a+b$. It
implies that $p=2$ and $d(u_1)=d(s_{p-1})=3$. Similarly, we get
$q=2$ and $d(u_n)=3$. By \eqref{s21} and \eqref{sab},
$S(u_1)=3+d(w_2)+d(u_2)=8$. Note that $d(w_2),d(u_2)=2,3$ or $4$ by
Lemma \ref{3no}~(ii). It follows that $d(u_2)=2, d(w_2)=3$ or
$d(u_2)=3, d(w_2)=2$. We may suppose that $d(u_2)=2, d(w_2)=3$
without loss of generality. Then $k=3$ and $d(u_i)=2$ for $2\leq i
\leq n-1$ by Lemma \ref{p}~(i). Hence $G$ has no pendent vertex.
This is a contradiction.

If $l\geq3$, then $d(r_2)=2$ or $4$. If $d(r_2)=2$, then
$d(s_{p-1})=6-d(r_2)=4$. Thus $S(s_{p-1})=10$ by \eqref{sab}. This
is impossible since $S(s_{p-1})=4-2+3+d(s_{p-2})\leq9$ by
\eqref{s21}. If $d(r_2)=4$, then $S(r_2)=10$ by \eqref{sab}. This is
also impossible since $S(r_{2})=4-2+3+d(r_{3})\leq9$ by \eqref{s21}.

Hence $d(v_1)=d(v_m)=a+b$. Dually,  we have $d(u_1)=d(u_n)=a+b$.

Applying Lemma \ref{3no}~(iv) with $C=v_1v_2\dots v_mr_{l-1}\dots
r_2v_1$, we have, without loss of generality, $m\geq3$ and
$d(v_i)\neq a+b$ for some $2\leq i\leq m-1$. So $d(v_i)=2$ for all
$2\leq i\leq m-1$ by Lemmas \ref{3no}~(ii) and \ref{p}.

Applying \eqref{suv} with $(v,u)=(v_1,x)$ and $(v,u)=(v_2,x)$,
respectively. We have
\begin{equation}\label{de}
a=\frac{d(r_2)+d(s_{p-1})-1}{a+b-1}\ \mbox{and}\ a=d(v_3),
\end{equation}
respectively.

We claim that $m\geq4$. Otherwise, suppose that $m=3$. Then
$a=d(v_3)=a+b$ and
$d(r_2)+d(s_{p-1})=a(a+b-1)+1\geq3(a+b)-2>2(a+b)$, which is
impossible since $d(r_2),d(s_{p-1})=2$ or $a+b$ by Lemma
\ref{3no}~(ii). Hence $m\geq4$. Therefore $m=4$ by Lemma
\ref{p}~(ii). Thus $a=d(v_3)=2$. It follows from \eqref{de} that
$d(r_2)+d(s_{p-1})=2(a+b)-1$. By Lemma \ref{3no}~(iii),
$d(r_2),d(s_{p-1})=2$ or $a+b$.

If $d(r_2)=d(s_{p-1})=2$, then $a+b=2.5$, a contradiction.

If $d(r_2)=d(s_{p-1})=a+b$, then $2(a+b)=2(a+b)-1$, a contradiction.

If $d(r_2)=2,d(s_{p-1})=a+b$, then $a+b=3$. So $a=2,b=1$.

By Lemma \ref{p}~(i), $d(s_{i})=d(s_{p-1})=3$ for $2 \leq i \leq
p-1$.

For the vertex $u_1$, $S(u_1)=3+d(u_2)+d(w_2)=7$ by \eqref{s21} and
\eqref{sab}. By Lemma \ref{3no}~(ii), $d(u_2),d(w_2)=2$ or $3$. It
follows that $d(u_2)=d(w_2)=2$. So $S(u_2)=3+d(u_3)=5$. Thus
$d(u_3)=2$ and $n\geq4$. Hence $n=4$ by Lemma \ref{p}~(ii).

Similarly, we have $d(w_3)=d(r_2)=2$ and $k=l=4$.

By \eqref{s21} and \eqref{sab}, $S(u_4)=4+d(t_2)=7$. So $d(t_2)=3$.
Hence $d(t_i)=3$ for $1\leq i \leq q-1$ by Lemma \ref{p}~(i).

Therefore $G\in {\mathcal G}_5$ (see Fig.~4), where
$max\{k_1,k_2\}\geq1$. It is easy to see that any graph $G\in
{\mathcal G}_5$ is $2$-walk $(2,1)$-linear.

If $d(r_2)=a+b,d(s_{p-1})=2$, then with a similar argument of the
case $d(r_2)=2,d(s_{p-1})=a+b$, we have $G\in {\mathcal G}_6$ is
$2$-walk $(2,1)$-linear (see Fig.~4), where $max\{k_1,k_2\}\geq1$.

Up to now, we have complete the proof of the Lemma. $\square$

\begin{lemma}\label{t8} Let $G\in{\mathcal
G}$ with  $G_0\in {\mathcal T}_{8}$. Then $G\cong H_{i}$ for
$i=24,\dots,30$ (see Fig.~3) or $G\in{\mathcal G}_j$ for $j=7,8$
(see Fig.~4).
\end{lemma}
\noindent {\bf Proof.} Let $G_0=T_{8}$ (see Fig.~2), where
$n,m,k,l,p,q\geq2$. For convenience, we set $u_1=v_1=w_1$,
$u_n=t_1=s_p$, $v_m=r_1=t_q$ and $w_k=s_1=r_l$.

{\bf Case 1}.  There is a vertex $v\in\{w_2,v_{2},u_2\}$ with
$d_{G_0}(v)=2$ and $d(v)=a+b$. Without loss of generality, suppose
that $v=w_{2}$. Then $k\geq3$ and $d(w_i)=a+b$ for $2\leq i \leq
k-1$ by Lemma \ref{p}~(i). In particular, $,a\geq2,a+b\geq3$ by
Lemma \ref{3no}~(iii).

Applying Lemma \ref{suv} with $(v,u)=(w_2,x)$, we have $
a=\frac{d(w_1)+d(w_3)-2}{a+b-1}$. By Lemma \ref{3no}~(ii), we have
$d(w_1),d(w_3)=3$ or $a+b$.

If $d(w_1)=d(w_3)=3\neq a+b$, then $a+b\geq4$ and
$a=\frac{4}{a+b-1}<2$, a contradiction.

If $d(w_1)=a+b,d(w_3)=3$ or $d(w_3)=a+b,d(w_1)=3$ and $a+b\neq 3$,
then $a=\frac{a+b+1}{a+b-1}=1+\frac{2}{a+b-1}<2$, a contradiction.

If $d(w_1)=d(w_3)=a+b$, then $a=2$. We claim that $a+b=3$.
Otherwise, let $a+b>3$. For the vertex $w_1$,
$S(w_1)=2(a+b)+b=a+b-3+a+b+d(u_2)+d(v_2)$. Thus
$d(u_2)+d(v_2)=a+b+1$. By Lemma \ref{3no}~(ii), $d(u_2),d(v_2)=2,3$
or $a+b$. Note that $a+b>3$. We have $d(u_2)=2,d(v_2)=3$ or
$d(u_2)=3,d(v_2)=2$ or $d(u_2)=d(v_2)=3$. Without loss of
generality, we consider the following two cases:

If $d(u_2)=2,d(v_2)=3$, then $a+b=4$ and $m=2$. Note that $a=2$. We
have $b=2$. If $k\geq4$, then $d(w_k)=d(w_1)=4$ by Lemma
\ref{p}~(i). If $k=3$, then $d(w_k)=d(w_3)=4$. So $d(r_3)=2 $ or $4$
by Lemma \ref{3no}~(ii). On the other hand,
$S(v_2)=4+d(r_2)+d(t_{q-1})=8$ by \eqref{s21} and \eqref{sab}. So
$d(r_2)=d(t_{q-1})=2$. By \eqref{s21} and \eqref{sab},
$S(r_2)=6=3+d(r_3)$. Thus $d(r_3)=3$, a contradiction.

If $d(u_2)=d(v_2)=3$, then $a+b=5$ and $m=n=2$. Note that $a=2$. We
have $b=3$. By \eqref{s21} and \eqref{sab},
$S(v_2)=5+d(r_2)+d(t_{q-1})=9$. Thus $d(r_2)=d(t_{q-1})=2$. Hence
$S(r_2)=7=3+d(r_3)$. This is impossible since $d(r_3)\in\{2,3,5\}$.

Therefore $a+b=3$ and $a=2,b=1$.

By \eqref{s21} and \eqref{sab}, $S(w_1)=7=3+d(v_2)+d(u_2)$. So
$d(v_2)=d(u_2)=2$. Hence $S(v_2)=5=3+d(v_3)$. It follows that
$d(v_3)=2$ and $m\geq4$. Therefore $m=4$ by Lemma \ref{p}~(ii).

Similarly, we have
$n=l=p=4,d(u_2)=d(u_3)=d(r_2)=d(r_3)=d(s_2)=d(s_3)=2$ and $d(t_i)=3$
for $2\leq i\leq q-1$.

Therefore $G\in{\mathcal G}_8$ (see Fig.~4), where
$max\{k_1,k_2\}\geq1$. It is easy to see that any graph
$G\in{\mathcal G}_8$ is 2-walk $(2,1)$-linear.

{\bf Case 2}. For any vertex  $v\in\{w_2,v_{2},u_2\}$, we have
$d_{G_0}(v)=3$ or  $d(v)=2$.

If $G$ has no pendent vertex, then $n,m,k,l,p,q\leq4$ by Lemma
\ref{31}. If $n=m=k=l=p=q=2$, then $G$ is regular. It is well known
that a graph is regular if and only if it has exactly one main
eigenvalues. Thus $max\{n,m,k,p,q\}\geq3$. Without loss of
generality, suppose that $n\geq3$. We consider the following two
cases:

{\bf Subcase 1}.  $n=3$. Then $S(u_2)=6$ and $m,k,l,p,q=2$ or $3$ by
Lemma \ref{31}.

If $m=k=2$, then $S(u_1)=8=S(u_3)=2+d(t_2)+d(s_{p-1})$. So
$d(t_2)=d(s_{p-1})=3$ and hence $p=q=2$. Similarly,
$S(u_1)=8=S(v_2)=6+d(r_2)$. Thus $d(r_2)=2$ and $l=3$. Therefore
$G=H_{24}$ (see Fig.~3). By \eqref{suv}, $H_{24}$ is 2-walk
$(2,2)$-linear. With a similar argument, we have:

If $m=2,k=3$ or  $m=3,k=2$, then $G=H_{25}$ is 2-walk $(1,4)$-linear
(see Fig.~3).

If $m=k=3$, then $G=H_{26}$ is 2-walk $(0,6)$-linear (see Fig.~3).

{\bf Subcase 2}. $n=4$.  Then $m,k,l,p,q=2$ or $4$ by Lemma
\ref{31}. With a similar argument of Subcase 1, we have the
following cases:

If $m=k=2$, then $G=H_{27}$ is 2-walk $(3,-1)$-linear (see Fig.~3).

If $m=2,k=4$ or  $m=4,k=2$, then $G=H_{28}$ is 2-walk $(2,1)$-linear
(see Fig.~3).

If $m=k=4$, then $G=H_{29}$ is 2-walk $(1,3)$-linear (see Fig.~3).

If $G$ has at least one pendent vertex $x$, then $x\in N(v)$ for
$v\in\{w_1,v_{m},u_{n},w_{k},r_{i_1},t_{i_2},s_{i_3}\}$, where
$2\leq i_1\leq l-1,2\leq i_2\leq q-1,2\leq i_3\leq p-1$. If
$v\in\{r_{i_1},t_{i_2},s_{i_3}\}$, then with a similar argument of
Case 1, we have $G\in\mathcal{G}_8$. Hence we may suppose that
$v\in\{w_1,v_{m},u_{n}\}$. In particular, assume that $x\in N(w_1)$
without loss of generality. Then $d(w_1)=a+b\geq4$ and
$d(u)=d_{G_0}(u)$ for $u\in\{r_{i_1},t_{i_2},s_{i_3}\}$, where
$2\leq i_1\leq l-1,2\leq i_2\leq q-1,2\leq i_3\leq p-1$. Applying
Lemma \eqref{suv} with $(v,u)=(w_1,x)$, we have\[
a=\frac{a+b-3+d(u_2)+d(v_2)+d(w_2)-(a+b)}{a+b-1}=
\frac{d(u_2)+d(v_2)+d(w_2)-3}{a+b-1}.\] Note that
$d(u_2),d(v_2),d(w_2)=2,3$ or $a+b$ and $a+b\geq4$. We consider the
following cases without loss of generality:

If $d(u_2)=d(v_2)=d(w_2)=2$, then $a=\frac{3}{a+b-1}<2$, a
contradiction.

If $d(u_2)=d(v_2)=2, d(w_2)=3$, then $a=\frac{4}{a+b-1}<2$, a
contradiction.

If $d(u_2)=d(v_2)=2, d(w_2)=a+b$, then $a=1+\frac{2}{a+b-1}<2$, a
contradiction.

If $d(u_2)=2,d(v_2)= d(w_2)=3$, then $a=\frac{5}{a+b-1}<2$, a
contradiction.

If $d(u_2)=2,d(v_2)= 3,d(w_2)=a+b$, then $m=2$ and
$a=1+\frac{3}{a+b-1}$. By Lemmas \ref{21} and \ref{3no}~(ii), we
have $a+b=4$ and $a=2$. Hence $d(w_i)=d(w_2)=4$ for $2\leq i\leq
k-1$. By \eqref{s21} and \eqref{sab}, $S(w_{k-1})=4+2+d(w_k)=10$.
Thus $d(w_k)=4$. Similarly, $S(v_2)=4+d(r_{2})+d(t_{q-1})=8$. It
implies that $d(r_{2})=d(t_{q-1})=2$. Thus $S(r_2)=6=3+d(r_3)$.
Hence $d(r_3)=3$. On the other hand, $d(r_3)\in\{2,4\}$ by Lemma
\ref{3no}~(ii) and the fact that $d(w_k)=4$. This is a
contradiction.

If $d(u_2)=2, d(v_2)=d(w_2)=a+b$, then $a=2+\frac{1}{a+b-1}$ is not
an integer, a contradiction.

If $d(u_2)=d(v_2)=d(w_2)=3$, then $n=m=k=2$ and $a=\frac{6}{a+b-1}$.
It follows from Lemma \ref{3no}~(iii) that $a+b=4$ and $a=b=2$. By
\eqref{s21} and \eqref{sab}, $S(v_2)=4+d(r_2)+d(t_{q-1})=8$. It
implies that $d(r_2)=d(t_{q-1})=2$. Hence $S(r_2)=3+d(r_3)=6$. It
follows that $d(r_3)=3$ and $l=3$. Similarly, we have $p=q=3$.
Therefore $G=H_{30}$ (see Fig.~3). By \eqref{suv}, $H_{30}$ is
2-walk $(2,2)$-linear.

If $d(u_2)=d(v_2)= 3,d(w_2)=a+b$, then $n=m=2$ and
$a=1+\frac{4}{a+b-1}$. Thus $a+b=5$ and $a=2$ by Lemma
\ref{3no}~(iii). Hence $d(w_1)=d(w_2)=a+b=5$. For the vertex $v_2$,
$S(v_2)=5+d(r_2)+d(t_{q-1})=9$. It implies that
$d(r_2)=d(t_{q-1})=2$. Thus $S(r_2)=3+d(r_3)=7$, which is impossible
since $d(r_3)\in\{2,3,5\}$ by Lemma \ref{3no}~(ii).

If $d(u_2)=3,d(v_2)= d(w_2)=a+b$, then $a=2+\frac{2}{a+b-1}$ is not
an integer, a contradiction.

If $d(u_2)=d(v_2)= d(w_2)=a+b$, then $n=m=k=2$ and $a=3$. We claim
that $p=l=2$. Otherwise, let $p,l>2$. By \eqref{s21},
$S(w_1)=a+b-3+3(a+b),S(w_2)=a+b-3+a+b+d(s_2)+d(r_{l-1})$. Note that
$d(s_2),d(r_{l-1})=2$ by assumption. We have $S(w_1)>S(w_2)$. On the
other hand, $S(w_1)=S(w_2)$ by \eqref{sab}, a contradiction. Hence
$p=l=2$. Similarly, we have $q=2$.  Thus $G\in {\mathcal G}_7$ (see
Fig.~4), where $b\geq1$. It is easy to see that any graph $G\in
{\mathcal G}_7$ is 2-walk $(3,b)$-linear.

Up to now, we have complete the proof of the Lemma. $\square$

\begin{theorem} The graphs $H_i$ for $i=1,\dots,30$ and those in ${\mathcal
G}_j$ for $j=1,\dots,8$ are all connected tricyclic graphs with
exactly two main eigenvalues.
\end{theorem}
\noindent {\bf Proof.}  It follows directly from Lemmas
\ref{t1}--\ref{t8}. $\square$

\end{document}